# *How to Conquer the Liar and Enthrone the Logical Concept of Truth: an Informal Exposition*

BORIS ČULINA
*University of Applied Sciences Velika Gorica, Velika Gorica, Croatia*

*This article informally presents a solution to the paradoxes of truth and shows how the solution solves classical paradoxes (such as the original Liar) as well as the paradoxes that were invented as counterarguments for various proposed solutions ("the revenge of the Liar"). This solution complements the classical procedure of determining the truth values of sentences by its own failure and, when the procedure fails, through an appropriate semantic shift allows us to express the failure in a classical two-valued language. Formally speaking, the solution is a language with one meaning of symbols and two valuations of the truth values of sentences. The primary valuation is a classical valuation that is partial in the presence of the truth predicate. It enables us to determine the classical truth value of a sentence or leads to the failure of that determination. The language with the primary valuation is precisely the largest intrinsic fixed point of the strong Kleene three-valued semantics (LIFPSK3). The semantic shift that allows us to express the failure of the primary valuation is precisely the classical closure of LIFPSK3: it extends LIFPSK3 to a classical language in parts where LIFPSK3 is undetermined. Thus, this article provides an argumentation, which has not been present in contemporary debates so far, for the choice of LIFPSK3 and its classical closure as the right model for the truth predicate. In the end, an erroneous critique of Kripke-Feferman axiomatic theory of truth, which is present in contemporary literature, is pointed out.*

## 1. *Introduction*

The concept of truth has various aspects and is a frequent subject of philosophical discussions. Philosophical theories usually consider the concept of truth from a wider perspective. They are concerned with questions such as—Is there any connection between the truth and the world? And, if there is—What is the nature of the connection?[1] Contrary to these theories, the analysis of the paradoxes of truth is of a logical nature because it deals with the internal semantic structure of a language, the mutual semantic connection of sentences, above all the connection of the sentences that speak about the truth of other sentences and the sentences whose truth they speak about. That is why every solution to the paradoxes of truth necessarily establishes a certain logical concept of truth.

The paradoxes of truth are "symptoms of disease" (Tarski 1969: 66): they show that there is a problem in our basic understanding of language, and they are a test for any proposed solution. Thereby, it is important to make a distinction between the *normative* and the *analytic* aspect of the solution.[2] The former tries to ensure that paradoxes will not emerge. The latter attempts to explain why paradoxes arise and to construct a solution based on that explanation. Of course, the practical aspect of the solution is also important. It tries to ensure a good framework for logical foundations of knowledge, for related problems in artificial intelligence and for the analysis of the natural language.

In the twentieth century, two solutions stood out, Tarski's (Tarski 1933, Tarski 1944) and Kripke's (Kripke 1975) solution. They initiated a whole series of considerations, from elaboration and critique of their solutions to proposals for different solutions. For the solution that is informally presented in this article, only Tarski's and Kripke's solutions are important, so other solutions will not be considered.[3]

Tarski's analysis emphasised the T-scheme as the basic intuitive principle for the truth predicate. According to Tarski, to examine the truth value of the sentence *"'snow is white' is a true sentence"*, we must examine whether snow is white. Thus, for the truth predicate the following must hold:

"snow is white" is a true sentence if and only if snow is white

This should be true for every declarative sentence $S$:

$\overline{S}$ is a true sentence if and only if $S$

where $\overline{S}$ is the name of the sentence $S$. For a particular sentence, we can always achieve this with quotation marks, as shown in the example of the sentence "snow is white". Tarski called this sentence scheme the

---

[1] A good overview of philosophical theories of truth can be found in (Glanzberg 2018). The author's position is set out in (Čulina 2020).

[2] In (Chihara 1979: 590), Chihara writes about "the preventative problem of the paradox" and about "the diagnostic problem of the paradox".

[3] An overview of various solutions can be found in (Beall et al. 2020).

*T-scheme*. However, if we apply the T-scheme to the sentence *L*: "$\overline{L}$ is a false sentence" (the famous Liar sentence), we will get a contradiction (the Liar paradox):

$\overline{L}$ is a true sentence if and only if $\overline{L}$ is a false sentence

Thus, Tarski's analysis showed the inconsistency of the T-scheme with the classical logic for the languages in which the Liar can be expressed, such as natural language.

Tarski's solution is to preserve the classical logic and to restrict the T-scheme to parts of the language. Tarski showed that if a language *L* meets some minimum requirements, we can consistently talk about the truth values of sentences of *L* only inside another "essentially richer" (Tarski's term) metalanguage *ML*. In *ML*, the T-scheme can only be set for the language *L*. This solution is in harmony with the idea of reflexivity of thinking and it has become very fertile for mathematics and science in general. For example, in chemistry, using the sentences of a language *L* we describe chemical processes, and using the sentences of *ML* we talk about the truth values of sentences of the language *L*.

Tarski does not deal with the analysis of the mechanism that leads to the paradoxes of truth, but only with the logical analysis of the formal inference of the contradiction. Not wanting to give up classical logic, the solution necessarily leads him to separate the metalanguage in which the T-scheme is expressed and the language for which the T-scheme is expressed, as a formal means of eliminating contradiction. Although Tarski does not explicitly say it anywhere, his solution suggests that the paradoxes of truth have their source in the violation of the reflexivity of thinking: talk about truth is an act of reflection whereby we transcend the original language. However, Tarski's solution is primarily of a normative nature. The mechanism of the paradoxes of truth is not analysed but paradoxes are blocked by a syntactic restriction. In *ML* we can speak only of the truth values of the sentences of the language *L*, so in *ML* the paradoxes of truth cannot be expressed at all. As for the liar paradox, the maximum approximation allowed by the syntactic restriction is the Limited Liar: when *L* is part of *ML*, under certain conditions, we can construct in *ML* the sentence

*LL*: $\overline{LL}$ is a false sentence of the language *L*[4]

If *LL* belonged to the language *L*, we could apply the T-scheme to *LL*:

$\overline{LL}$ is a true sentence of the language *L* if and only *LL*

According to the construction of the sentence *LL*, we get a contradiction:

$\overline{LL}$ is a true sentence of the language *L* if and only $\overline{LL}$ is a false sentence of the language *L*

---

[4] For example, this can be realized if for ML we choose the language of Peano's arithmetic and for L we choose $\Sigma_n$ sentences of the language (Kaye 1991: 126).

However, from this contradiction it follows that *LL* does not belong to the language *L*. Thus, *LL* is certainly not a false sentence of the language *L*. So, it is a false sentence of the language *ML*.

Kripke showed that there is no natural syntactic restriction to the T-scheme as set out in Tarski's solution, but that we must look for the solution in the semantic structure of language. Consider the first example given by Kripke (Kripke 1975: 690). In the New Testament Saint Paul writes:

> One of Crete's own prophets has said it: "Cretans are always liars, evil brutes, idle bellies". He has surely told the truth.

In accordance with Tarski's approach, we can take as an object language the language composed of all the declarative sentences uttered by the Cretans together with the above statements of Saint Paul. In doing so, we will consider Saint Paul's first sentence to be true, which is an acceptable assumption. We will also assume that Saint Paul said all the above. For the sake of simpler expression, the sentence "Cretans are always liars, evil brutes, idle bellies" will be called "What one of Crete's own prophets said". The second St Paul's sentence is context dependent, so we will explicate it as the sentence "What one of Crete's own prophets said is true" and call it "What Saint Paul said". The application of the T-scheme for the object language gives us here:

1) What one of Crete's own prophets said is true if and only if Cretans are always liars, evil brutes, idle bellies
2) What Saint Paul said is true if and only if What one of Crete's own prophets said is true

According to 1), if What one of Crete's own prophets said is true then Cretans are always liars. So, What one of Crete's own prophets said is a lie. From this contradiction we conclude that What one of Crete's own prophets said is not true. By 2), we further conclude that What Saint Paul said is not true either. There is nothing paradoxical in the analysis so far.[5] However, let us consider what we can infer from the fact that What one of Crete's own prophets said is not true. By 1), it follows that Cretans are not always liars, evil brutes, idle bellies. So, we learned something about Cretans. It may seem odd that we have concluded something factual based on the T-scheme. However, we used Saint Paul's first statement as a factual assumption about the Cretans, and the T-scheme was only part of the logical mechanism by which we derived the above factual statement about the Cretans from this assumption. From a logical point of view, everything seems to be fine. However, we can imagine the extreme situation: that "one of Crete's own prophets" is the only Cretan, that he is not an evil brute or idle belly. That would mean he sometimes tells the truth. But we can go further and imagine that he made only one claim in his entire life—the one Saint Paul mentions. That would mean that What one of Crete's own prophets said is a true statement. And so, we got a contradiction again. In such a situation we are given a paradox: What one of Crete's own prophets said is true if and only if it is false, and so What Saint Paul said is true if and only if it is false.

---

[5] Except perhaps for those who believe that everything written in the New Testament must be true.

In his article, Kripke describes a much more realistic situation in which the statements made have a certain truth value in normal conditions, but under some specific conditions they become paradoxical. In Kripke's words (Kripke 1975: 691):

> many, probably most, of our ordinary assertions about truth and falsity are liable, if the empirical facts are extremely unfavourable, to exhibit paradoxical features.

Kripke's analysis clearly showed that for a language in which one sentence speaks about the truth values of other sentences, what is expected and what is paradoxical in the language cannot be separated on the syntactic or internal semantic level: it depends on the reality that the language is talking about, and not on the way we use the language. Thus, according to Kripke, it is necessary to include this risk in the theory of truth. Sentences that speak of the truth values of other sentences, although syntactically correct and meaningful, under some conditions depending on the reality to which the language refers may not make a determinate claim about that reality: they will not give a classical truth value, *True* or *False*. Then we assign the third value to them: *Undetermined*. The meaning of the third value is simply that the sentence has no classical truth value. Such an analysis leads to the study of languages with partial two-valued semantics, which, by introducing *Undetermined* as the third value, is technically equivalent to the study of languages with three-valued semantics.

Kripke did not give any definite model. He gave a theoretical framework for investigations of various models—each fixed point in each monotone three-valued semantics can be a model for the truth predicate. Each such model gives a natural restriction on the T-scheme: the T-scheme is valid for all sentences that have a classical truth value in that model, while for the others it is undetermined.[6] However, as with Tarski, the proposed solutions are normative—we can express the paradoxical sentences, but we escape a contradiction by declaring them undetermined.

Kripke took some steps in the direction of finding an analytic solution. He preferred the strong Kleene three-valued semantics (SK3 semantics below) for which he wrote it was "appropriate" but did not explain why it was appropriate. One reason for such a choice is probably that Kripke finds paradoxical sentences meaningful. This elimi-nates the weak Kleene three-valued semantics which in the standard interpretation[7] corresponds to the idea that paradoxical sentences are meaningless, and thus undetermined. Another reason could be that the SK3 semantics has the so-called investigative interpretation. According to this interpretation, this semantics corresponds to the classical determination of truth values, whereby all sentences that do not have an already determined value are temporarily considered undetermined. When we determine the truth values of some of these sentences, then we can also determine the truth values of some of the sentences that are composed of them, which were undetermined until then. For example, if we know that *S* is a true sentence and we do not yet know the truth of the sentence *T*, then according to the classical truth valuation of the conjunction, we do not yet know the truth of the sentence *S and T* (we will know it only when we know the truth value of the sentence *T*) but we do know that the disjunction *S or T* is true. This truth valuation corresponds exactly to the SK3 semantics.[8] Kripke supplemented this investigative interpretation with an intuition about learning the concept of truth in the presence of the truth predicate. That intuition deals with how we can teach someone who is a competent user of an initial language (without the truth predicate "to be true") to use sentences that contain the truth predicate. That person knows which sentences of the initial language are true and which are not. We give her the rule to assign the attribute "to be true" to the former and deny that attribute to the latter. In that way, some new sentences that contain the truth predicate, and which were undetermined until then, become determined. So, the person gets a new set of true and false sentences with which she continues the procedure. This intuition leads directly to the minimal fixed point of the SK3 semantics (MIFPSK3 below) as an analytically acceptable model for the concept of truth.

---

[6] For Kripke, as well as for my further analysis, the rules associated with the T-scheme are much more important than the T-scheme itself: that whenever the sentence S has a truth value, then the sentence "is true" has the same value and vice versa.

In the structure of the fixed points of a language with the truth predicate, two fixed points stand out, the minimal fixed point and the largest intrinsic fixed point. The first has the structural property that every sentence that has a classical truth value at the minimal fixed point has the same value at other fixed points. The largest intrinsic fixed point has the structural property that it is the largest fixed point such that every sentence that has a classical truth value in it has no opposite classical value at any other fixed point (it is compatible with all other fixed points). Kripke's work gives an internal characterisation of MIFPSK3, which follows from Kripke's description of the learning process of the concept of truth: at that fixed point, only those sentences whose truthfulness is based on the described learning process have a truth value. Starting with Kripke, the largest intrinsic fixed point is mostly mentioned as an interesting solution because of its structural properties. Kripke writes (Kripke 1975: 709):

> The largest intrinsic fixed point is the unique "largest" interpretation of T($x$) which is consistent with our intuitive idea of truth and makes no arbitrary choices in truth value assignments. It is thus an object of special theoretical interest as a model.

---

[7] Some philosophers have given a different interpretation of the weak Kleene three-valued semantics, e.g. (Beall 2016).

[8] In weak Kleene three-valued semantics, if T would be a meaningless sentence, there is no need for further truth valuation, because automatically all sentences containing T are also meaningless.

Since then, nothing much has changed in philosophical debates. Thus, forty years later, Horsten in his review article (Horsten 2015) writes:

> Until now, the intrinsic fixed points have not been investigated as intensively as they should perhaps be.

In (Čulina 2001) and in PhD thesis (Čulina 2004) I gave an analytic solution to the problem of the paradoxes of truth. In (Čulina 2001) it has been shown that this solution is precisely the largest intrinsic fixed point of the SK3 semantics (LIFPSK3 below) together with its classical closure. In this way, LIFPSK3 got a specific interpretation. This article provides an argumentation, which has not been present in contemporary philosophical discussions, for the choice of LIFPSK3 and its classical closure as the right model for the logical concept of truth. The solution will be informally described, and it will be demonstrated how it solves the classical paradoxes of truth (such as the original Liar) as well as the paradoxes that have been invented as counterarguments for various solutions to the paradoxes of truth ("the revenge of the Liar"). I will try to make the argumentation as simple as possible, so that the consideration can be followed by someone who does not have any special knowledge of the techniques related to Tarski's and Kripke's analysis. Finally, one of the confirmations of the naturalness of the solution of the problem of the logical concept of truth should be that such a solution can be explained in simple language, understood and used by any interested language user who does not have a special mathematical and philosophical education. All these informal considerations can be formalised by the means developed in (Čulina 2001). Some parts of the text are taken from (Čulina 2001) and PhD thesis (Čulina 2004). (Čulina 2001) contains the formal results to which this argumentation refers, while the PHD thesis contains the basic elements of the argumentation itself. However, much of what is only stated there has been elaborated and supplemented here to present rounded and convincing argumentation for the logical concept of truth introduced in these works.

## 2. *An analysis of the paradoxes of truth*

An analysis of the paradoxes of truth will be done on sentences. Tarski and Kripke state the technical reasons for this choice. In (Tarski 1944: 342) Tarski writes:

> By "sentence" we understand here what is usually meant in grammar by "declarative sentence"; as regards the term "proposition", its meaning is no-toriously a subject of lengthy disputations by various philosophers and logicians, and it seems never to have been made quite clear and unambiguous. For several reasons it appears most convenient to apply the term "true" to sentences, and we shall follow this course.

Kripke writes (Kripke 1975: 691):

> I have chosen to take sentences as the primary truth bearers not because I think that the objection that truth is primarily a property of propositions (or "statement") is irrelevant to serious work on truth or to the semantic paradoxes. On the contrary, I think that ultimately a careful treatment of the problem may well need to separate the "expresses" aspect (relating sentences to propositions) from the "truth" aspect (putatively applying to propositions). ... The main reason I apply the truth predicate directly to linguistic objects is that for such objects a mathematical theory of self-reference has been developed.

A convincing argument for choosing sentences for truth bearers was given by Quine in (Quine 1986: 1). This choice has an undoubted technical advantage because the subject of study is specific language forms, and not abstract objects of unclear nature. It is also a reflection of my deep conviction that language is not just a means of writing down and communicating thoughts but an essential part of thinking.[9]

Roughly, by "classical language" will be meant every language which is modelled upon the everyday language of declarative sentences. Due to definiteness, a language of the first order logic, which has an explicit and precise description of form and meaning, will be considered. By "language" will be meant an interpreted language, a language form together with an interpretation.

The interpretation of a first order language determines the external semantic structure of the language, a connection between the language and the subject matter of the language. The connection is based on external assumptions on the language use: (i) the language has its own domain of interpretation—a collection of objects that the language speaks of, (ii) every constant denotes an object, and every variable in a given valuation denotes an object, (iii) every function symbol symbolises a function which applied to objects gives an object, (iv) every predicate symbol symbolises a predicate which applied to objects gives a truth value, *True* or *False*. For simplicity, I will assume that the language has names for all objects in its domain. In doing so, $\bar{a}$ will be the name for an object $a$.

The inner organisation of a first-order language is determined by the rules of the construction of more complex language forms from simpler ones, starting with names, variables and function symbols for building terms, and with atomic sentences for building sentences. In these constructions we use special symbols which identify the type of the construction. With each construction, and thus the symbol of the construction, a semantic rule is associated that determines the semantic value of the constructed whole using the semantic values of the parts of the construction.[10] These rules determine the internal semantic structure of the language. The symbol of a language construction will be termed *logical symbol* or *logical constant* if the associated semantic rule is an internal language rule: the rule does not refer to the reality the language speaks of, except possibly referring to external assumptions of the language use. For example, connectives and quantifiers are logical symbols of the language.[11] I will further argue that the truth predicate "to be true" is also a logical symbol of the language. The interconnectedness of the truth values of sentences of a language belongs to one aspect of the concept of truth which I will term the *logical aspect of the concept of truth*. The connection of the truth value of a sentence with reality forms, for example, the second aspect of the concept of truth.[12] Since the paradoxes of truth occur in the context of the logical aspect of the concept of truth, I believe that each of their solutions establishes a certain logical concept of truth.

---

[9] My view of the essential role of language in thinking and rational cognition is explained in (Čulina 2021a).

The external assumptions of the language use have grown from everyday use of language where we are accustomed to their fulfilment, but there are situations when they are not fulfilled. The Liar paradox and other paradoxes of truth are witnesses of such situations for the external assumption (iv). Let's consider the sentence $L$ (the Liar):

$L$: $\overline{L}$ is a false sentence. (or "This sentence is false.")

Using the usual understanding of language, to investigate the truth value of $L$ we must investigate what it says. But it says precisely about its own truth value, and in a contradictory way. If we assume it is true, then it is true what it says—that it is false. But if we assume it is false, then it is false what it says, that it is false, so it is true. Therefore, it is a self-contradictory sentence. What is disturbing is the paradoxical situation that we cannot determine its truth value.

The same paradoxicality, but without contradiction, emerges in the investigation of the following sentence $I$ (the Truth-teller):

$I$: $\overline{I}$ is a true sentence. (or "This sentence is true.")

Contrary to the Liar to which we cannot associate any truth value, to this sentence we can associate the truth as well as the falsehood with equal mistrust. There are no additional specifications which would make a choice between the two possibilities.

I will begin the analysis of the paradoxes of truth with a basic observation that the above paradoxical sentences are meaningful because we understand well what they say, even more, we used that in the unsuccessful determination of their truth values. However, they witness the failure of the classical procedure for the truth value determination in some "extreme" situations. According to the classical procedure, the examination of the truth value of a sentence is reduced to the examination of the truth values of the sentences from which it is constructed according to the classical truth value conditions for that type of construction. Thus, for example, the examination of the truth value of a sentence of the form φ *or* ψ is reduced to the examination of the truth values of the sentences $\varphi$ and $\psi$. The reduction is performed according to the truth value conditions for the logical connective *or*: $\varphi$ or $\psi$ is true when at least one of the sentences $\varphi$ and $\psi$ is true, and false when both $\varphi$ and $\psi$ are false sentences. Likewise, a sentence of the form $\forall x\, P(x)$ (where $\forall$ is the standard symbol for *for all*) is a true when the sentences $P(\overline{a})$ are true for every object $a$ from the domain of the language, and it is false when $P(\overline{a})$ is false for at least one object $a$. Thus, the examination of the truth value of a sentence comes down to the examination of the truth value of the sentences from which it is constructed (if these sentences contain free variables, then we must look at all valuations of these variables). Examining the truth values of these sentences is in the same way reduced to examining the truth values of the sentences from which they are constructed, etc.

---

[10] In a given interpretation and a given valuation of variables, the semantic value of a term is the object described by the term and the semantic value of a sentence is its truth value.

[11] In (Čulina 2021b) the concept of logical symbol of a language is elaborated in more detail.

[12] In (Čulina 2020) various aspects of the concept of truth are analysed.

We can visualise this procedure on the graph whose nodes are sentences of the language, where each sentence points with an arrow to the sentences to which, according to the classical truth value conditions, the examination of its truth value is reduced. Each type of sentence construction gives the corresponding type of elementary block of such a graph. To illustrate, the blocks corresponding to the constructions using negation (*not*), the disjunction (*or*), and the universal quantor ($\forall$) are shown below:

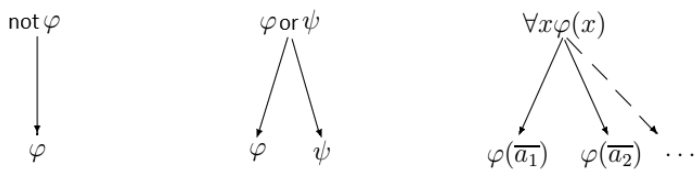

Each sentence has its own *semantic graph* to which the sentence is a distinguished node, and the graph is composed of all sentences on which, according to the truth value conditions, the truth value of a given sentence hereditarily depends.[13]

To determine the truth value of a given sentence, according to the classical truth value conditions, we must investigate the truth values of all sentences to which it points, then possibly, for the same reasons, the truth values of the sentences to which these sentences point, and so on. Every such path along the arrows of the graph leads to atomic sen-

---

[13] The semantic graph of the whole language can be defined analogously. The semantic graphs of individual sentences are its subgraphs.

tences (because the complexity of sentences decreases along the path). In situations where a language doesn't talk about the truth values of its own sentences, the truth values of its atomic sentences don't depend on the truth values of some other sentences. The atomic sentences are the leaves of the semantic graph of the given sentence. To investigate their truth values, we must investigate external reality they are talking about. The classical assumption of a language is that every atomic sentence has a definite truth value. So, the procedure of determination of the truth value of the given sentence stops and gives a definite truth value, *True* or *False*. Formally, this is secured by the recursion principle which says that there is a unique function from sentences to truth values, which obeys the classical truth value conditions and its values on atomic sentences are identical to externally given truth values.[14] Such is, for example, the language of a scientific field, but not the everyday language in which there are frequent discussions about the truthfulness of claims made by others. In such situations, the above analysis can be, and is, disrupted when atomic sentences use the truth predicate to speak of the truth values of other sentences of the language. These are sentences of the form $T(\overline{\varphi})$, where "T" is the symbol for the truth predicate "to be true", and $\overline{\varphi}$ is the name of a sentence $\varphi$ of the language. Such an atomic sentence is not a leaf of a semantic graph, but points with an arrow to the sentence $\varphi$ on which its truth value depends:

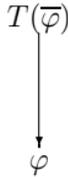

The truth value conditions of this construction are the basic conditions of the logical concept of truth: that $T(\overline{\varphi})$ is true when $\varphi$ is true, and $T(\overline{\varphi})$ is false when $\varphi$ is false. Given that the semantic rule of this construction is an internal semantic rule (it connects the truths of sentences $T(\overline{\varphi})$ and $\varphi$ independently of the reality the language speaks of), the truth predicate is a logical symbol of the language, in the same way that, for example, connectives and quantifiers are logical symbols of the language. In this sense, it is perfectly correct to speak of this concept of truth as a *logical concept of truth*. The only difference in relation to connectives and quantifiers is in universality. Only a language that has its own sentences in the domain of its interpretation (possibly through coding) can have a logical symbol of its own truth predicate.

---

[14] Note that, even when we know the true values of the leaves, this procedure is generally not computable because although the semantic graph of a given sentence has finite depth (the reduction to the leaves takes place in the finite number of steps), the leaves themselves can be infinitely many.

In the presence of the truth predicate, it can happen that the procedure of determination of the truth value of a given sentence does not stop at atomic sentences but, under the truth value conditions of the truth predicate, continues through each atomic sentence of the form T($\overline{\varphi}$) to the sentence $\varphi$. Because of the possible "circulations" or other kinds of infinite paths, there is nothing to ensure the success of the procedure. Truth paradoxes just witness such situations. Five illustrative examples follow.

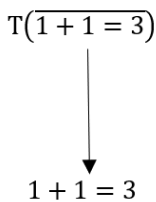

The procedure of the truth value determination has stopped on the atomic sentence for which we know is false, so T($\overline{1+1=3}$) is false, too.

The Liar: For *L*: T($\overline{not\ L}$) we have

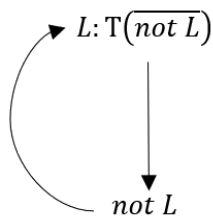

But now the procedure of the truth value determination has failed because the truth value conditions can't be fulfilled. The truth value of T($\overline{not\ L}$) depends on the truth value of *not L* and this again on *L*: T($\overline{not\ L}$) in a way which is impossible to obey.

The Truth-teller: For *I*: T($\overline{I}$) we have

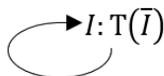

Now, there are, as we have already seen, two possible assignments of the truth values to the sentence *I*. But this multiple fulfilment we must consider as a failure of the classical procedure, too, because the procedure assumes to establish a unique truth value for every sentence.

The Logician: *Log*: T($\overline{Log}$) or T($\overline{not\ Log}$) (This sentence is true or false)

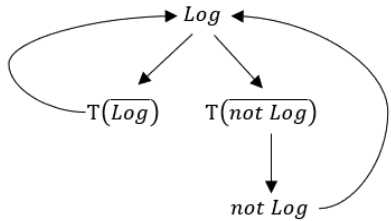

If *Log* were false then, by the truth conditions, T($\overline{not\ Log}$) would be false, not *Log* would be false too, and finally *Log* would be true. Therefore, such valuation of the graph is impossible. But if we assume that *Log* is true, the truth conditions generate a unique consistent valuation. Therefore, the truth determination procedure gives the unique answer—that *Log* is true.

The law of excluded middle for the Truth-teller: *I or not I*

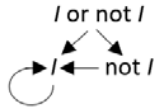

Now there are two truth valuations of the semantic graph of the sentence *I or not I*. In both valuations, it takes the same value: *True*. However, in one valuation the sentence *I* takes the value *True* and in the other *False*. Given that the classical procedure requires that not only the initial sentence, but every sentence included in the examination, if it has a truth value, then has a unique truth value, we must also consider this situation as a failure of the classical procedure for determining the truth value of the sentence *I or not I*. Having failed to determine the truth value of sentence *I*, we have not been able to determine the truth value of the sentence *not I*, and therefore neither of the sentence *I or not I*.

The paradoxes of truth emerge from a confrontation of the implicit assumption of the success of the classical procedure of the truth value determination and the discovery of the failure. As previous examples show such assumption is an unjustified generalisation from common situations to all situations. We can preserve the classical procedure, also the internal semantic structure of the language, but we must reject universality of the assumption of its success. The awareness of that transforms paradoxes to normal situations inherent to the classical procedure. I consider this the diagnosis of the paradoxes of truth.

## 3. *The proposed solution*

The previous diagnosis shows us the way to the solution—the formulation of the *partial two-valued semantics* of language which, when the procedure of determining the truth value of a given sentence gives a unique truth value, *True* or *False*, attaches that value to the sentence, and when the procedure fails, it does not attach any truth value to the sentence. This kind of semantics can be described as the *three-valued semantics* of language—simply the failure of the procedure will be declared as the third value (*Undetermined*). It has not any additional philosophical charge. It is only a convenient technical tool for the description. In formulating the partial two-valued semantics, we will start from these properties:

1) The semantics coincides with the classical semantics on atomic sentences whose truth values are determined by the external reality they are talking about.
2) In the semantics all sentences are meaningful.
3) The semantics has classical truth value conditions for connectives and quantifiers.
4) In the semantics $T(\overline{\varphi})$ is true when $\varphi$ is true, and false when $\varphi$ is false (a variant of the T-scheme).
5) When the classical procedure of determining the truth value of a given sentence assigns it a unique truth value then the semantics assigns that value to the sentence, otherwise it does not assign a truth value to the sentence.

Properties 1) and 4) need no comment. Property 2) was commented at the beginning of this analysis. The fact that we cannot determine the truth values of paradoxical sentences does not mean that they are not meaningful. We understand their meaning quite well. Moreover, we use this meaning essentially in the (unsuccessful) determination of their truth values. The consequence of this property is that all sentences have meaning, regardless of whether some part of the sentence is paradoxical or not. Otherwise, as soon as one part of the sentence was paradoxical, the whole sentence would be meaningless.[15] Here is one argument as to why it is not an acceptable solution to consider paradoxical sentences to be meaningless. If we were to accept that paradoxical sentences have no meaning, it would make no sense to determine their truth values. Thus, we could not determine which sentences are paradoxical, i.e., they have no meaning.[16]

For property 3) it is only important to note that the rejection of the success of the classical procedure of the truth value determination doesn't change the meaning of the classical truth value conditions. They are stated in a way independent of the assumption that sentences must have a truth value. They specify the truth value of a compound sentence in terms of the truth values of its direct components regardless of whether they have truth values or not. The lack of some truth value may lead, but does not have to, to the lack of the truth value of the compound sentence. For example, the truth value conditions of the sentence $\varphi$ *and* $\psi$ are: $\varphi$ *and* $\psi$ is true when both $\varphi$ and $\psi$ are true, and false when at least one of the sentences $\varphi$ and $\psi$ is false. It says nothing about the existence of the truth values of $\varphi$ and $\psi$, but only sets conditions among the truth values. The functioning of the truth value conditions in the new situation is illustrated by the example of the following sentences (where $L$ is the Liar, and $I$ is the Truth-teller):

$$L \text{ or } 0 = 0$$

---

[15] This would lead to the weak Kleene three-valued semantics of the language.

[16] Thus this argument rejects the weak Kleene three-valued semantics as a solution to the paradoxes of truth.

By the classical truth value conditions for the connective *or*, this sentence is true precisely when at least one of the basic sentences is true. Because 0=0 is true consequently the total sentence is true, although $L$ has not a truth value. Equally, if we apply the truth value conditions on the connective *and* to the sentence

$$L \text{ and } 0 = 0$$

the truth value will not be determined. Namely, for the sentence to be true both basic sentences must be true, and it is not fulfilled. For it to be false at least one basic sentence must be false, and this also is not fulfilled. So, non-existence of the truth value for $L$ leads to non-existence of the truth value for the whole sentence. Let's analyse

$$I \text{ or not } I$$

Since $I$ does not have a truth value, *not I* does not have a truth value, so *I or not I* also does not have a truth value.

Property 5) is a key property. It expresses the basic idea of this approach: the lesson of the paradoxes of truth is that the classical procedure of determining the truth value does not have to succeed. By failure we mean that, respecting the classical conditions of truth, we cannot assign a truth value to a sentence, or we can assign two truth values to it. However, property 5) stated in this way is not precise enough because the classical determination of truth values is not an algorithmic process and in concrete situations we manage to implement it in various ways. Furthermore, rejecting the assumption of the existence of a unique truth value complicates the process, because now it is not enough to find one valuation of a given sentence, but it is necessary to examine whether there are other valuations, not only of the given sentence but also of other sentences included in the examination. That's why we must give property 5) a more objective formulation that does not talk about the real or idealised process of determining the truth values of sentences, but about the existence of these values. In (Čulina 2001) it is shown that:

*There is a unique partial two-valued semantics (association of truth values to sentences) with the following properties*:
1) The semantics coincides with the classical semantics on atomic sentences whose truth values are determined by the external reality they are talking about.
2) The semantics has classical truth value conditions for connectives and quantifiers.
3) In the semantics T($\overline{\varphi}$) is true when $\varphi$ is true, and false when is false (a variant of the T-scheme).
4) The semantics is unique in the sense that on the set of all sentences to which it associates truth values, every other semantics that fulfils the previous three conditions does not associate different truth values (it can happen that the other semantics does not associate a truth value to some of these sentences).
5) The semantics is the largest such semantics in the sense that on the set of sentences to which it does not assign truth values, every other semantics that fulfils the previous four conditions also does not assign truth values.

Below, I will call this semantics *the partial two-valued semantics*. This is exactly the requested extension of classical semantics to situations where it is not guaranteed that every sentence is true or false because this semantics accurately identifies when the classical procedure of determining the truth value of a sentence will succeed and when it will not.

This result gives a "license" to the classic procedure of determining the truth value of a sentence in situations where not all sentences have a truth value. The restriction of the partial two-valued semantics to the semantic graph of a given sentence is the truth valuation of the graph, which the classical procedure should determine by its success or failure. Thereby, the classical procedure does not need to determine the entire valuation, but only that part that is sufficient to determine the truth value of a given sentence or to determine that it has no truth value. For example, to determine the truth value of the sentence $\exists x\, \varphi(x)$ (where $\exists$ is the standard symbol for exists), if among all sentences of the form $\varphi(\overline{a})$ we find one that is true then we do not have to examine the others, nor do we have to worry about whether any of them is undetermined. Likewise, when we know that some sentences are undetermined, we can use this in determining non-existence of the truth values of other sentences. For example, for the sentence *L and* 0=0, knowing that *L* is undetermined allowed us to conclude that *L and* 0=0 is also undetermined. Thereby, not only the truth value conditions for the connective *and* do not give us the truth value for *L and* 0=0 but the failure of the classical procedure in determining the truth value of *L* leads to the failure of determining the truth value of *L and* 0=0. This example shows that not only the classical truth value conditions of the conjunction of two sentences do not depend on whether these sentences have a truth value but the conditions also determine how the failure of the determination of truth values is propagated. It is easy to see that this is also true in general: all the truth value conditions not only determine the connection between truth values but also determine how the failure of the determination is propagated. If we look at the associated three-valued semantics, it is not difficult to show that these are precisely the conditions of the SK3 semantics. Thus, SK3 have a special interpretation here: the SK3 conditions are the classical truth value conditions supplemented by the conditions of propagation of the failure to determine truth values.

In (Čulina 2001) it was proved that in the labyrinth of literature on the paradoxes of truth (Beall et al. 2020), the partial two-valued semantics described above is positioned as the largest intrinsic fixed point of the SK3 semantics (LIFPSK3) with a specific interpretation. In that way, the above presented argumentation for the partial two-valued semantics is also the argumentation for the choice of LIFPSK3 among all fixed points of all monotone three-valued semantics for the right model of the logical concept of truth.

In (Kremer 1988: 245), Kremer writes:

> Within Kripke's theoretical framework there are two leading candidates for the "correct" interpretation of the truth predicate: the minimal fixed point and the largest intrinsic fixed point. ... We are thus led to distinguish two plausible versions of the principle of the supervenience of semantics. First, there is the view of that the correct interpretation of truth is the minimal fixed point; as we saw, this has often been taken to be "Kripke's theory of truth". Second, there is the view that the largest intrinsic fixed point is the correct interpretation of truth. Unfortunately for the champion of supervenience, there seem to be considerations in support of both of these views.

I will give some arguments as to why I consider LIFPSK3 with the interpretation described in this article to be a better solution than MIFPSK3 with Kripke's interpretation. The main argument concerns the content-wise interpretations of these fixed points. In Kripke, it is an interpretation of learning the concept of truth in the presence of the truth predicate, here an interpretation of determining truth values of sentences that language users actually do.

In Kripke, the SK3 semantics has an investigative interpretation: while we have not yet determined the truth values of some sentences, they are undetermined. In the process of learning the concept of truth in the presence of the truth predicate, more and more sentences gain truth value. So, some hitherto undetermined sentences become determined, which, according to the truth value conditions, entails that some others sentences become determined. However, some sentences will remain undetermined forever. Thus, as Visser noted in (Visser 1989: 651), the SK3 interpretation changes: "not yet" interpretation of undetermined value in the learning process (because we haven't learned the concept of truth enough yet), in MIFPSK3 becomes "not ever" interpretation (a sentence is undetermined because its truth val-ue can never be determined through the process of learning the concept of truth). In the interpretation developed in this article, undetermined sentences are those sentences to which the classical procedure of determining truth values does not give a unique truth value. SK3 naturally derives from the classical procedure of determining truth values, which in the presence of the truth predicate is not always successful. In this interpretation, SK3 is simply the classical semantics complemented by the propagation of its own failure. For example, let's analyse the Liar in both interpretations. In Kripke's interpretation, learning the predicate of truth, we will not face the Liar at any level to determine its truth value. The Liar is simply inaccessible to us in that process and, if we strictly adhere to the metaphor of learning the predicate of truth, we will never even know that the Liar is inaccessible to us. In contrast, the interpretation developed in this article provides only a formal framework for the process a language user actually undertakes when encountering the Liar and examining its truth value. The language user will easily determine that the Liar is undetermined.

Furthermore, in Kripke's interpretation, language users learn the truth predicate in an extensional way, collecting more and more sentences that fall under the predicate and more and more sentences that do not fall under the predicate. However, as explained in this article, the truth predicate is a logical concept: it is determined by the internal semantics of language and it should not be learned experientially, just as, for example, the logical meaning of the connective *and* should not be learned experientially. As we know the meaning of the connective *and* when we are given its truth value conditions, so we know the meaning of the truth predicate, when we are given its truth value conditions: $T(\overline{\varphi})$ is true when $\varphi$ is true, and it is false when $\varphi$ is false. From this definition of the logical concept of the truth predicate, which is a variant of the T-scheme, and which corresponds to the basic intuition of the language user, arises the interpretation developed in this article which gives LIFPSK3. In Kripke's case the opposite is true: from the intuition about learning the truth predicate follows MIFPSK3 as an extensional a posteriori definition of the truth predicate.

Finally, to learn someone which sentences to associate with the predicate "is true", we must first know it ourselves. According to Kripke's interpretation, someone should first learn us which sentences to associate with the predicate "is true". Thus, the idea of learning the truth predicate leads to an infinite regress: learning grounded truth is ungrounded.

That the aspect of learning the concept of truth and understanding the concept of truth is not one and the same, Yablo has already noted in (Yablo 1982: 118), but in the context of MIFPSK3:

> If the inheritance aspect is the one lying behind the attempt to picture grounding in terms of the learning of 'true', then the dependence aspect is the one behind the attempt to picture grounding in terms of the understanding of 'true'. What do we do when we have to evaluate a sentence—say "The sentence 'Snow isn't white' is true" or the sentence "The sentence 'Snow is white' is true' is not true"—involving complicated attributions of truth? Evidently, we try to figure out what its truth-value depends on, and then what that depends on, and so on and so forth in the hope of eventually making our way down to sentences not containing 'true' which can be evaluated by conventional means. … But the fact that the majority of those who grappled with grounding before Kripke tended to see things from the standpoint of dependence suggests that there is something intuitively satisfying about the dependence approach.

As already commented, MIFPSK3 and LIFPSK3 have distinguished structural properties in the structure of all fixed points. Kripke's description of the learning process gave a characterisation of MIFPSK3 independent of other fixed points. The analysis developed in this article provides a characterisation of LIFPSK3 that is also independent of other fixed points. However, while Kripke's characterisation is global—the learning process yields all the truths and falsehoods of MIFPSK3—the LIFPSK3 characterisation developed here is local: the truth value determination of a given sentence takes place only on the semantic graph of the sentence. The characterisation of LIFPSK3 developed here, and not the Kripke's characterisation of MIFPSK3, corresponds to the way a language user determines the truth value of a sentence. Starting from a given sentence, the language user tries to determine its truth value by examining its semantic graph, and not by collecting more and more true and false sentences according to the instructions for learning, and hoping that the given sentence will appear in one of those groups. In this letter case, as it was illustrated above on the example of the Liar, the language user can never determine that a sentence is undetermined: it is constantly in the "not yet" interpretation and can never switch to the "not ever" interpretation.

LIFPSK3 contains MIFPSK3 as a subset, which can also be considered an advantage of LIFPSK3. To all sentences, which MIFPSK3 assigns a truth value, LIFPSK3 also assigns this same value. But in addition, LIFPSK3 assigns truth values to sentences that are undefined in MIFPSK3. Such, for example, is the sentence the Logician, which is assigned the value *True* by the classical truth determination procedure, as shown above, while it is undetermined in MIFPSK3. Of course, a remark can be made here that MIFPSK3 is a better choice for this very reason, because in MIFPSK3 truth values are given only to those sentences whose truth value is "grounded" in the reality that the language speaks about. Such is not, for example, the Logician, but the sentence *L or* $0 = 0$ is: this sentence is true because the atomic sentence $0 = 0$ is true. Although Kripke formally calls grounded all sentences that have a truth value in MIFPSK3, on an intuitive level they are grounded because their truth values are determined by the process of learning the concept of truth which starts from the truth values of atomic sentences that speak of reality. However, as already stated, paradoxes of truth fall under the logical concept of truth that connect the truth values of sentences independently of the reality that the language speaks of, so the solution to the paradoxes should not include reality. The classical procedure for determining truth values does not require us to examine whether a sentence is grounded or not, but only whether we can associate a unique truth value with it or not. Of course, in this examination, we can arrive at atomic sentences that talk about external reality, but we don't have to, as the example of the Logician shows. In this comparison it is also seen that the choice between Kripke's interpretation of MIFPSK3 and the interpretation of LIFPSK3 described in this paper is a choice between non-logical and logical concept of truth.

The next section will show that some of Gupta's critiques (Gupta 1982) of fixed points apply to MIFPSK3 but not to LIFPSK3. Thus, the critiques turn into the argument that LIFPSK3 is a more acceptable model for the truth predicate than MIFPSK3.

So, for now we have two semantics of a language with the truth predicate. We have the *classical* or *naive semantics* in which paradoxes occur because this semantics assumes that each sentence is true or false, i.e., it assumes that the process of determining truth values always gives an unambiguous answer. And we have its repair to the two-valued partial semantics of the language, i.e., to the three-valued semantics of the language, which accepts the possibility of failure of the classical procedure of determining truth values. I will call this semantics the *primary semantics* of the language. However, to remain on the partial two-valued semantics would mean that the logic would not be classical, the one we are accustomed to. Concerning the truth predicate itself, it would imply the preservation of its classical logical sense in the two-valued part of the language extended by the "silence" in the part where the classical procedure fails. For example, the T-scheme is true only for sentences that have a classical truth value. For other sentences it is undetermined. Although in a meta-description, $T(\overline{\varphi})$ has the same truth value (in the three-valued semantic frame) as $\varphi$, that semantics is no longer the initial classical semantics (although it extends it) nor it can be expressed in the language itself: the language is silent about the third value. Or better said, the third value is the reflection in a meta-language of the silence in the language. So, the expressive power of the language is weak. For example, the Liar is undetermined. Although we have easily said it in the metalanguage, we cannot express in the language itself, because, as it has already been said (in the metalanguage), the Liar is undetermined. Not only that this "zone of silence" is unsatisfactory for the above reasons (it leads to the three-valued logic, it loses the primary sense of the truth predicate and it weakens the expressive power of the language), but it can be overcome by a natural *additional valuation* of the sentences which emerges from recognising the failure of the classical procedure. "Natural", in the sense that it is precisely this move that a language user makes in the end when faced with the failure of the classical procedure. This point will be illustrated on the example of the Liar. On the intuitive level of thinking, by recognising the Liar is not true nor false we state that it is undetermined. However, this is not a claim of the original language but of the metalanguage in which we describe what happened in the language. Moreover, in the metalanguage, we can continue to think. Since the Liar is undetermined, it is not true what it claims—that it is false. Therefore, the Liar is false. But this does not lead to restoring of the contradiction because by moving to the metalanguage we have made a *semantic shift* from the primary partial two-valued semantics (or the three-valued semantics) toward its two-valued meta description. Namely, the Liar talks of its own truth in the frame of the primary semantics, while the last valuation is in the frame of another semantics, which I will term the *final semantics* of the language. The falsehood of the Liar in the final semantics doesn't mean that it is true what it says (that it is false) because the semantic frame is not the same. The falsehood of the Liar in the final semantics means that it is false (in the final semantics) what the Liar talks of its own primary semantics: that it is false in the primary semantics (because it is undetermined in the primary semantics). So, not only have we gained a contradiction in the naive semantics, i.e., the third value in the primary semantics, but we also have gained additional information about the Liar.

A key element in the above consideration is a semantic shift in thinking. It is closely related to the reflexivity of thinking, which appears in two variants in the paradoxes of truth. The first variant takes place in the primary semantics of the language, and the second in the transition to the final semantics of the language. In the first variant, the reflexivity of thinking occurs in the transition from the use of the sentence $S$ to the mention of the sentence $S$. The most significant example of this transition in the context of the paradoxes of truth is the transition from the statement $S$ to the statement "$\overline{S}$ is a true sentence". Thereby, two aspects of the concept of truth should be distinguished. To examine whether "snow is white" is a true sentence, we must investigate reality, see what colour the snow is. So, the truth value of that sentence depends on reality. Therefore this aspect of the concept of truth is not of a logical nature. To examine the truth value of the sentence ""snow is white" is a true sentence" we must examine the truth value of the sentence "snow is white". Thus, the truth value of that sentence also depends (indirectly) on reality. However, the truth predicate only articulates this transfer of truth, just as, for example, the connective *and* articulates the transfer of the truth values of a conjunction to the conjuncts. According to the classical truth conditions on the truth predicate, the predicate connects the truth values of two sentences in the primary semantics in a way that is independent of the reality the language is talking about. That is why the truth predicate is a logical symbol in the primary semantics and falls under the logical notion of truth in the primary semantics.

Another variant of reflexivity of thinking occurs at the level of the whole language—in the transition from the original language to the metalanguage by which we describe the original language. In the context of the paradoxes of truth, this transition was illustrated above on the example of the Liar, when we concluded that the Liar is undetermined in the primary semantics. This conclusion belongs to the metalanguage by which we describe truth valuations in the original language. The metalanguage has the same syntax as the original language (we will see below that the predicate "is undetermined" can be defined by the predicate "is true"), but not the same semantics: it has a different connection between the truth values of sentences, and this connection is a classical two-valued semantics, called here the final semantics of the language. The metalanguage is a classical language with classical semantics and an external reality that it talks about in the same way that a classical language is, for example, the language we use to describe car engines. Only here it is not so obvious, because the external reality that the metalanguage talks about is another language that has the same syntax as the metalanguage (but not the same semantics). The key element of the semantic distinction between these two languages is the truth predicate that we have analysed so far. Given that we now have two languages, we must first express this predicate of truth more precisely, to make it clear that it is the truth predicate of the original language. That is why instead of "is true" we will now use "is true in the primary semantics". This does not change its role in primary semantics one bit—the logical role described above. In the final semantics, the sentence "$\overline{S}$ is true in the primary semantics" has the same meaning as in the primary semantics: it asserts that the sentence $S$ is true in the primary semantics. However, in the final semantics the truth predicate of the primary semantics connects the truth value of the sentence $S$ in the primary semantics with the truth value of the sentence "$\overline{S}$ is true in the primary semantics" in the final semantics. The semantic mechanism here is the same as with the predicate "is a diesel motor", which connects the engine type $x$ with the truth value of the statement "$x$ is a diesel motor". However, since in the final semantics the truth predicate connects the truth values of two semantics, and not engines and truth values, confusion can easily occur if we don't take care which truth value belongs to which semantics. Just as the predicate "is a diesel motor" is not a logical symbol of language, because we must examine the external reality of language—a given engine—to determine the truth of the corresponding sentence, so too, the truth predicate "is true in the primary semantics" is not a logical symbol of the metalanguage (the final semantics) because we have to investigate the external reality of the metalanguage—investigate the truth value of a sentence in the original language (the primary semantics)—to determine the truth value of that sentence in the final semantics. The semantic shift that allowed us to complete the analysis of the Liar paradox, and that allows us to complete the analysis of other paradoxes of truth, as will be shown later, is precisely this change of the role of the truth predicate of the primary semantics from the logical symbol of the primary semantics, which connects the truth values within the primary semantics, into a non-logical symbol of the final semantics, which connects the truth values of the primary and the final semantics. This change leads to a change in the overall semantics of the language—from the partial two-valued primary semantics to the classical two-valued final semantics of the language.

It is easy to legalise this intuition about semantic shift. Sentences of a language with the truth predicate will always have the same meaning, but the language will have two valuation schemes—the primary and the final truth valuation. In both semantics the meaning of the truth predicate is the same: T($\overline{\varphi}$) means that $\varphi$ is true in the primary semantics. But the valuation of the truth value of the atomic sentence T($\overline{\varphi}$) is different. While in the primary semantics the truth value conditions for T($\overline{\varphi}$) are classical (the truth of T($\overline{\varphi}$) means the truth of $\varphi$, the falsehood of T($\overline{\varphi}$) means the falsehood of $\varphi$, and consequently T($\overline{\varphi}$) is undetermined just when $\varphi$ is undetermined), in the final semantics it is not so. In it, the truth of T($\overline{\varphi}$) means that $\varphi$ is true in the primary semantics, and falsehood of T($\overline{\varphi}$) means that $\varphi$ is not true in the primary semantics. It does not mean that it is false in the primary semantics, but that it is false or undetermined. So, formally looking, in the final semantics T($\overline{\varphi}$) inherits truth from the primary semantics, while other values transform to falsehood. That is why we say that this semantics is the *classical semantic closure* of the primary semantics, or in full terminology, the classical semantic closure of LIFPSK3. Due to the monotonicity of the primary semantics this means that the final semantics supplements the primary semantics in the area of its silence. If a sentence in the primary semantics has a classical value (*True* or *False*), it will have that value in the final semantics as well. If a sentence is undetermined in the primary semantics (a paradoxical sentence) then it will have a classical truth value in the final semantics that just carries information about its indeterminacy in the primary semantics. Therefore, the final semantics is the classical two-valued semantics of the language that has for its subject precisely the primary semantics of the language, and it extends the primary semantics in the part where the primary semantics is silent, using just the information about the silence.

We can see best that this is a right and a complete description of the valuation in the primary semantics by introducing predicates for other truth values in the primary valuation:

> F($\overline{\varphi}$) ("$\varphi$ is false in the primary semantics") ↔ T($\overline{not\ \varphi}$)
> U($\overline{\varphi}$)("$\varphi$ is undetermined in the primary semantics") ↔ *not* T($\overline{\varphi}$) *and not* F($\overline{\varphi}$)

According to the truth value of the sentence $\varphi$ in the primary semantics we determine which of the previous sentences are true and which are false in the final semantics. For example, if $\varphi$ is false in the primary semantics then $F(\overline{\varphi})$ is true while others ($T(\overline{\varphi})$ and $U(\overline{\varphi})$) are false in the final semantics. Once the final two-valued valuations of atomic sentences are determined in this way, the final valuation of every sentence is determined by means of the classical truth value conditions and the principle of recursion. This valuation not only preserves the primary logical meaning of the truth predicate (as the truth predicate of the primary semantics) but it also coincides with the primary valuation where it is determined.

I think that when a language user is confronted with a paradox of truth, his thinking ends in this final semantics. Therefore, the solution to the paradoxes of truth should include this semantics. Although both the primary and the final semantics share the same linguistic forms, the final semantics is the minimum metalanguage for the primary semantics by which we complete the analysis of paradoxical situations.

In the end of his article, Kripke warns that the complete description of paradoxical situations in a language with the truth predicate belongs to a metalanguage which has its own concept of truth, so the analysis of the concept of truth with fixed points remains incomplete, as in Tarski's approach. Kripke writes (Kripke 1975: 714):

> The necessity to ascend to a metalanguage may be one of the weaknesses of the present theory. The ghost of the Tarski hierarchy is still with us.

I do not think that the existence of a metalanguage with its concept of truth means that the analysis conducted here is incomplete. As already discussed, such a view arises from mixing various aspects of the concept of truth. The aim of this analysis is the logical concept of truth described on the page 11. It differs from the aspect of the concept of truth that is most important to us—truth that discriminates what is and what is not in the world that a language speaks of. The latter aspect of truth belongs to the external semantics of the language, its connection with the world, while this logical aspect of the concept of truth belongs to the internal semantics of the language. The critique of resorting to a metalanguage cannot be applied to the logical concept of truth because the truth values we associate with sentences of the metalanguage do not fall under the logical concept of truth. In particular, the concept of truth in the final semantics is not a logical concept of truth. It is equal to the concept of truth in other sciences. Of course, as in the languages of mechanical engineering, the question of the truth of sentences in the final semantics can be discussed in an appropriate metalanguage (and I've been doing it all along in these considerations). But this is a different type of problem than the problem of paradoxical sentences.[17]

---

[17] This is a problem of the truth regress: whenever we express a statement, we express its truth value with another statement whose truth value we express with another statement, etc.

## 4. *Conquering the Liar*

Having in mind this double semantics of the language (triple, if we also count the classical naive semantics), we can easily solve all truth paradoxes. On an intuitive level we have already done it for the Liar:

$L$: $F(\overline{L})$ ("This sentence is false.")

The form of the solution is always the same. A paradox in classical thinking means that the truth value of a sentence is undetermined in the primary semantics. But then it becomes an information in the final semantics with which we can conclude the truth value of the sentence in the final semantics. To make it easier to track solutions to other paradoxes, I will sometimes distinguish by appropriate prefixes what the truth valuation is about: I will put prefix "p" for the primary semantics and prefix "f" for the final semantics. In that way we will distinguish for example "f-falsehood" and "p-falsehood".

The Strengthened Liar is "the revenge of the Liar" for solutions that seek a way out in truth value gaps, i.e., in the introduction of the third value—*Undetermined*:

$SL$: *not* $T(\overline{SL})$ ("This sentence is not true.")

In the classical semantics it leads to a contradiction in the same way as the Liar because there "not to be true" is the same as "to be false". The paradox is used as an argument against the third value in the following way (e.g., in (Burge 1979)). If we accept that The Strengthened Liar takes on the value *Undetermined*, it means that what it is saying is true—that it is not true (but undetermined)—and so the contradiction is renewed. However, the last step is wrong because a semantic shift has occurred. The conclusion that The Strengthened Liar is undetermined is the conclusion in the final semantics. So, when we say in the end that what he says is true, this is the concept of truth of the final semantics, while the concept of truth The Strengthened Liar mentions is the concept of truth of the primary semantics. So, the truth of the final semantics is that The Strengthened Liar is not true in the primary semantics.

It is interesting that the whole argumentation can be done directly in the final semantics, not indirectly by stating the failure of the classical procedure. The argumentation is the following. If $SL$ were f-false, then it would be f-false what it said—that it is not p-true. So, it would be p-true. But it means (because the final semantics extends the primary one) that it would be f-true and it is a contradiction with the assumption. So, $SL$ is f-true. This statement does not lead to a contradiction but to an additional information. Namely, it follows that what $SL$ talks about is f-true—that it is not p-true. So, it is p-false or p-undetermined. If it were p-false it would be f-false too, and this is a contradiction. So, it is p-undetermined.

Note that, although the Liar and the Strengthened Liar are both p-undetermined, the latter is f-true while the former is f-false.

In (Burge 1979), Burge introduces the following the revenge of the Liar for the truth value gaps solutions:

$BL$: $F(\overline{BL})$ or $U(\overline{BL})$ ("This sentence is false or undetermined.")

When we consider it in the classical semantics, if it were true then it would be false or undetermined, which is a contradiction. If it were false, then it would be true—again a contradiction. So, again we make a semantic shift and in the final semantics we conclude that it is undetermined. This means that in the final semantics it is true. Or, if we express ourselves with prefixes, that sentence is p-undetermined and f-true.

The semantic shift in argumentation is best seen in the following variant, the so-called Metaliar:

1. The sentence on line 1 is not true.
2. The sentence on line 1 is not true.

The sentence on line 1 is The Strengthened Liar so it is undetermined. If we understand the second sentence as reflection on the first sentence, which we have determined to be undetermined, then the second sentence is true. So, it turns out that one and the same sentence is both undetermined and true. In (Gaifman 1992), Gaifman uses this example to motivate the association of truth values not with sentences as sentence types but with sentences as sentence tokens. Thus, Gaifman solves the paradox by separating the same sentence type into two tokens of which the first is undetermined and the second true. In my approach, it is precisely the separation of the primary and the final semantics of the same sentence. In the 1st line it gets the undetermined value in the primary semantics, while in the 2nd, by reflection on the primary semantics, it gets the value *True* in the final semantics.

In (Skyrms 1984), Skyrms introduced the Intensional Liar, to point out the intensional character of the Liar. Namely, if in The Strengthened Liar

(1): (1) is not true.

we replace (1) with the standard name of the sentence denoted by that sign, we get the sentence

"(1) is not true" is not true.

While sentence (1) is undetermined, this harmless substitution seems to have given us the sentence which is not undetermined but true (because "(1) is not true" is undetermined, and so it is not true). But here, too, there has been a semantic shift in the truth valuation that we can clarify with prefixes:

" "(1) is not true" is not p-true" is f-true.

## 5. *Conquering the companions of the Liar*

In the same way, paradoxes that have a different type of failure of the classical procedure, such as the Yablo paradox (Yablo 1993), are solved. Consider the following infinite set of sentences ($i$), $i \in \mathbb{N}$:

($i$) For all $k > i$ ($k$) is not true.

If the sentence ($i$) were true, then all the following sentences would not be true. But that would mean on the one hand that ($i+1$) is not true, and on the other hand, since all the sentences after it are not true, that ($i+1$) is true. So, all the above sentences are not true. But if we look what they claim entails that they are all true. This contradiction in the classical semantics turns into a true claim of the final semantics that all these sentences are p-undetermined. From what they say about their primary semantics, as with the Strengthened Liar, it follows that they are all f-true.

That the solution of the problem of the paradoxes of truth presented here is not related to negation will be illustrated by the example of Curry's paradox (Curry 1942):

$C$: $T(\overline{C}) \to l$ ("If this sentence is true then $l$")

where $l$ is any false statement. On the intuitive level, if $C$ were false then its antecedent $T(\overline{C})$ is false, and so the whole conditional $C$ is true: we got a contradiction. If $C$ was true then the whole conditional ($C$) and its antecedent $T(\overline{C})$ would be true, and so the consequent $l$ would be true, which is impossible with the choice of $l$ as a false sentence. Therefore, we conclude in the final semantics that $C$ is p-undetermined, and so it is f-true (because the antecedent is f-false).

All the paradoxical sentences analysed above led to contradictions in the classical semantics. Thus, in the final semantics, we concluded that they are undetermined in the primary semantics, from which we further determined their truth value in the final semantics. We could also analyse them directly in the final semantics, as was done with the Strengthened Liar. There, the contradiction would turn into a positive classical two-valued argumentation by which we would determine its truth value in both the primary and the final semantics. However, the situation is different with paradoxes which do not lead to a contradiction, which permit more valuations, like the Truth-teller. The analyses of the Truth-teller gives that it is p-undetermined. It implies that it is not p-true which means that ($I$: $T(\overline{I})$) it is not f-true. So, $I$ is f-false. However, although the conclusion is formulated in the final semantics, the reasoning that led to that conclusion cannot be formulated in the final semantics because it involves the analysis of the corresponding semantic graph. Of course, if we enrich the metalanguage with the description of semantic graphs and their truth value valuations then we could translate the whole intuitive argumentation into the final semantics.

In (Gupta 1982), Gupta gave several arguments against Kripke's fixed points. The solution presented here includes LIFPSK3, so this critique also applies to the solution developed in this article.

One of Gupta's criticisms, which has already been present in the literature, is that not all classical laws of logic are valid in fixed points. E.g., for a language containing the Liar, the logical law $\forall x \, not \, (T(x) \, and \, not \, T(x))$ is undetermined in each fixed point of the SK3 semantics (if we choose the Liar for $x$, we get the undetermined sentence). But since the analysis of paradoxes cannot avoid the presence of sentences that have no classical truth value, the analysis naturally leads to a three-valued language for which we cannot expect the logical laws of a two-valued language to apply. However, the SK3 semantics is maximally adapted to the two-valued logic: the logical truths of the two-valued logic are always true in SK3 when they are determined. Furthermore, the transition to the final semantics definitely solves this problem because that semantics is two-valued, and $\forall x \, not \, (T(x) \, and \, not \, T(x))$ is true in this semantics.

A somewhat more inconvenient situation is that $\forall x \, not \, (T(x) \, and \, not \, T(x))$, like other logical laws, is not true in the minimal fixed point even when there is not the Liar like or the Truth-teller like sentences. Namely, then the stated logical law is not true for its own sake—to determine its truth, the truth of all sentences, including itself, must be examined. In this way it can be seen that it is an ungrounded sentence, i.e., undetermined in MIFPSK3. But in LIFPSK3, it is true. We can easily check this by trying to give it a classic truth value. Namely, to examine its truth, we must examine whether the condition $not \, (T(x) \, and \, not \, T(x))$ is valid for each sentence $x$. Since we assume that language has no paradoxical sentences, it is only necessary to examine whether this is true of the law itself. If the law is false, then this condition is true of the law, so the law itself is true: we get a contradiction. Thus, the law must be true, and it is easy to show that this truth value does not lead to contradiction. Since the procedure of determining a truth value has assigned a unique truth value to this logical law, it is true in LIFPSK3. It means that this Gupta's critique turns into an argument for LIFPSK3.

The second type of Gupta's critique seeks to show that some quite intuitive considerations about the concept of truth are inconsistent with the fixed points of SK3 semantics. Gupta constructed the following example in (Gupta 1982) (Gupta's paradox). Let us have the following statements of persons $A$ i $B$:

$A$ says:

(a1) Two plus two is three. (false)
(a2) Snow is always black. (false)
(a3) Everything $B$ says is true. ( )
(a4) Ten is a prime number. (false)
(a5) Something $B$ says is not true. ( )

*B* says:

>    (b1) One plus one is two. (true)
>    (b2) My name is *B*. (true)
>    (b3) Snow is sometimes white. (true)
>    (b4) At most one thing *A* says is true. ( )

Sentences (a1), (a2), (a4), (b1), (b2) and (b3) are determined in each fixed point. However, (a3) and (a5) "wait" (b4), and (b4) "waits" them and so those sentences remain undetermined in the minimal fixed point. But on an intuitive level, it is quite easy for them to determine the classical truth value. Since (a3) and (a5) are contradictory, and all other statements of *A* are false, (b4) is true. But this means that (a3) is true and (a5) is false. However, this intuition coincides with the truth valuation in LIFPSK3. Thus, this Gupta's critique also turns into an argument for LIFPSK3. To find an intuitive counterexample for LIFPSK3 as well, Gupta replaces (a3) and (a5) with the following statements:

>    (a3\*): (a3\*) is true. ( )
>    (a5\*): "(a3\*) is not true" is true. ( )

Now at LIFPSK3, (a3\*) and (a5\*), and thus (b4), are undetermined. Gupta considers that on an intuitive level (b4) is true, because at most one of (a3\*) and (a5\*) is true. But in this step Gupta made a semantic shift from the primary semantics to the final, so (b4) is a true statement in the final semantics. This devalues his argument against LIPSK3.

## 6. *An erroneous critique of Kripke-Feferman theory*

In this last section I would like to draw attention to one erroneous critique of Kripke-Feferman axiomatic theory of truth (KF) which is present in contemporary literature, for example, in two contemporary respectable books on formal theories of truth. The models of this theory are the classical semantic closures of the fixed points of the SK3 semantics, and so the final semantics described in this paper, too.

In (Horsten 2011: 127) is the following text:

> So far, it seems that KF is an attractive theory of truth. However, we now turn to properties of KF that disqualify it from ever becoming our favourite theory of truth.
> Corollary 70: KF ⊢ $L \land \neg T(L)$, where *L* is the [strengthened] liar sentence.[18]
> …
> In other words, KF proves sentences that by its own lights are untrue. This does not look good. To prove sentences that by one's own lights are untrue seems a sure mark of philosophical unsoundness: It seems that KF falls prey to the strengthened liar problem.

---

[18] $\land$ and $\neg$ are the standard symbols for *and* and *not*.

In (Beall et al. 2018: 76) is the following text:

> But on the properties of truth itself, KF also has some features some have found undesirable. One example (discussed at length in Horsten 2011) is that KF ⊢ λ ∧ ¬Tλ. Unlike FS, KF gives us a verdict on Liars. But it seems to then deny its own accuracy, as it first proves λ, and then denies its truth. This makes the truth predicate of KF awkward in some important ways.

Both quoted texts repeat KF's critique dating back to Reinhardt (Reinhardt 1986), that axiomatic KF theory without additional restrictions is not an acceptable theory of truth. This means that its models, the classical semantic closures of the fixed points of SK3, are not acceptable solutions to the concept of truth. The reason is that the theory proves both The Strengthened Liar and that The Strengthened Liar is not true. The error in this reasoning stems from the indistinguishability of the primary (fixed point) and the final (classical semantic closure of the fixed point) semantics. KF has the role of axiomatically organising the final semantics, and what KF deduces are the true statements of the final semantics about the truth values of the primary semantics. We have already seen that The Strengthened Liar *SL* is true in the final semantics. Since KF axioms are valid in the final semantics, that KF ⊢ *SL* is not awkward but testifies to the strength of KF in the description of the fixed points. Furthermore, since *SL* is true in the final semantics, it means that it is not true in the primary semantics. So, that KF ⊢ ¬T($\overline{SL}$) is also not awkward but testifies to the strength of KF. These claims (in fact one claim KF ⊢ *SL* ∧¬T($\overline{SL}$)) are not contradictory, because different concepts of truth are involved.

## *References*

Beall, J. 2016. "Off-topic: A New Interpretation of Weak-Kleene Logic." *Australasian Journal of Logic* 13 (6).

Beall, J., Glanzberg, M., and Ripley, D. 2018. *Formal Theories of Truth*. Oxford: Oxford University Press.

Beall, J., Glanzberg, M., and Ripley, D. 2020. "Liar Paradox." In E. N. Zalta (ed.). *The Stanford Encyclopaedia of Philosophy*. Metaphysics Research Lab, Stanford University.

Burge, T. 1979. "Semantical Paradox." *Journal of Philosophy* 76: 169–198.

Chihara, C. 1979. "The Semantic Paradoxes: A Diagnostic Investigation." *Philosophical Review* 88 (4): 590–618.

Čulina, B. 2001. "The Concept of Truth". *Synthese* 126: 339–360.

Čulina, B. 2004. *Modelling the Concept of Truth Using the Largest Intrinsic Fixed Point of the Strong Kleene Three Valued Semantics* (*in Croatian Language*). PhD thesis. https://philpapers.org/archive/CULMTC.pdf.

Čulina, B. 2020. "The Synthetic Concept of Truth." Unpublished.

Čulina, B. 2021a. "The Language Essence of Rational Cognition with Some Philosophical Consequences". *Tesis* (*Lima*), 14 (19): 631–656.

Čulina, B. 2021b. "What is Logical in First Order Logic." Unpublished.

Curry, H. B. 1942. "The Inconsistency of Certain Formal Logics". *Journal of Symbolic Logic*, 7: 115–117.

B. Čulina, *How to Conquer the Liar* 31

Gaifman, H. 1992. "Pointers to Truth." *Journal of Philosophy* 89: 223–261.

Glanzberg, M. (ed.) 2018. *The Oxford Handbook of Truth*. Oxford: Oxford University Press.

Gupta, A. 1982. "Truth and Paradox." *Journal of Philosophical Logic* 11: 1–60.

Horsten, L. 2011. *The Tarskian Turn. Deflationism and Axiomatic Truth*. Cambridge: MIT Press.

Horsten, L. 2015. "One Hundred Years of Semantic Paradox." *Journal of Philosophical Logic* 44: 681–695.

Kaye, R. 1991. *Models of Peano Arithmetic*. Oxford: Clarendon Press.

Kremer, M. 1988. "Kripke and the Logic of Truth." *Journal of Philosophical Logic* 17 (3): 225–278.

Kripke, S. A. 1975. "Outline of a Theory of Truth." *Journal of Philosophy* 72: 690–716.

Quine, W. V. 1986. *Philosophy of Logic: Second Edition*. Cambridge: Harvard University Press.

Reinhardt, W. N. 1986. "Some Remarks on Extending and Interpreting Theories with a Partial Predicate for Truth." *Journal of Philosophical Logic* 15 (2): 219–251.

Skyrms, B. 1984. "Intensional Aspects of Semantical Self-Reference". In R. L. Martin (ed.). *Recent Essays on Truth and the Liar Paradox*. Oxford: Oxford University Press.

Tarski, A. 1933. "Pojęcie prawdy w językach nauk dedukcyjnych". *Towarzystwo Naukowe Warszawskie*. German translation, "Der Wahrheitsbegriff in den formalisierten Sprachen", Studia philosophica 1, 1935, 261–405.

Tarski, A. 1944. "The Semantic Conception of Truth". *Philosophy and Phenomenological Research* 4: 341–376.

Tarski, A. 1969. "Truth and Proof". *Scientific American* 220 (6): 63–77.

Visser, A. 1989. "Semantics and the Liar Paradox". In D. M. Gabbay et al. (eds.). *Handbook of Philosophical Logic*. Volume 4. Reidel: Springer.

Yablo, S. 1982. "Grounding, Dependence, and Paradox." *Journal of Philosophical Logic* 11 (1): 117–137.

Yablo, S. 1993. "Paradox without Self-Reference." *Analysis* 53: 251–252.